\documentclass[12pt]{article}
\usepackage{latexsym}
\usepackage{amssymb}
 
\newcommand{\proof}{{\noindent \bf Proof. }}

\newtheorem{defi}{Definition}
\newtheorem{lem}{Lemma}
\newtheorem{cor}{Corollary}
\newtheorem{prop}{Proposition}
\newtheorem{conj}{Conjecture}

\newcommand\mbf[1]{\mbox{\boldmath$#1$}}
\newcommand{\N}{{\mathbb N}}

\date{}

\begin{document}
\begin{titlepage}
\title{\bf GRAPH--DIFFERENT PERMUTATIONS}
{\author{{\bf  J\'anos K\"orner}
\\{\tt korner@di.uniroma1.it}
\\''La Sapienza'' University of Rome
\\ ITALY
\and {\bf Claudia Malvenuto}
\\{\tt claudia@di.uniroma1.it}
\\''La Sapienza'' University of Rome
\\ ITALY
\and {\bf G\'abor Simonyi}\thanks{Research partially supported by the
Hungarian Foundation for Scientific Research Grant (OTKA) Nos.\ 
T046376, AT048826, and NK62321.}
\\{\tt simonyi@renyi.hu}
\\R\'enyi Institute of Mathematics, Budapest
\\ HUNGARY}}

\maketitle
\begin{abstract}

We strengthen and put in a broader perspective previous results of the first
two authors on colliding permutations. 
The key to the present approach is a new non-asymptotic invariant for 
graphs. 
\end{abstract}
\end{titlepage}

\section{Introduction}

In \cite{KM} the first two authors began to investigate the following
mathematical puzzle. Call two permutations
of $[n]:=\{1,\dots,n\}$ {\em colliding} if, represented by linear orderings of
$[n]$, they 
put two consecutive elements of $[n]$ somewhere in the same
position. For the maximum cardinality $\rho(n)$ of a set of 
pairwise colliding permutations of $[n]$ the following conjecture was
formulated.  

\begin{conj}\label{conj:clak} {\rm (\cite{KM})}
For every $n\in \N$ 
$$\rho(n)={n \choose {\left \lfloor \frac{n}{2} \right \rfloor}}.$$
\end{conj} 

It was proved that the right hand side expression above is
actually an upper 
bound for $\rho(n)$, while the best lower bound given in \cite{KM} was a
somewhat deceiving  
\begin{equation}\label{eq:prev}
35^{{n/7}-O(1)}\leq \rho(n).
\end{equation}
 
The initial motivation for the present paper was to improve on the above lower
bound. For this purpose we will put the original problem in a broader
perspective leading to a new graph invariant that we believe to be interesting
on its own.  
For brevity's sake let us call a graph {\it natural} if its vertex set is a
finite subset of $\N$, the set of all positive integers and if the graph is
simple (without loops and multiple edges). An infinite permutation of $\N$ is
simply a linear ordering of all the elements of  
$\N$. 
(Instead of {\it infinite permutations of $\N$} we will often say simply
{\it infinite permutations} in the sequel.)
For an arbitrary natural graph $G=(V(G), E(G))$ we will 
call the
infinite permutations  
$\pi=(\pi(1), \pi(2), \dots, \pi(n), \dots)$ and $\sigma=(\sigma(1), \sigma(2),
\dots, \sigma(n), \dots)$ 
{\em $G$--different} if there is at least one $i\in \N$ for which 
$$\{\pi(i), \sigma(i)\}\in E(G).$$
(We will use the same expression for a
pair of finite sequences if at some coordinate they contain the two endpoints
of an edge of $G$.) 
Let $\kappa(G)$ be the maximum cardinality of a set of infinite permutations
any two elements of which are $G$--different. (It is easy to see that the
finiteness of $G$ implies that this number is finite as well, see
Lemma~\ref{lem:fin} below.)   
Clearly, the value of $\kappa$ is equal for isomorphic natural graphs.
In this paper we will analyze this quantity for some elementary graphs and
will apply some of the results to  
improve on the earlier estimates on $\rho(n)$. We have been able to determine
the value of $\kappa(G)$ only for some very small or simply structured graphs
$G$.  
Thus, to further simplify matters, we  
ask questions about the extremal values of $\kappa$ for graphs with a fixed
number of 
edges (and, eventually, vertices).  
We define

\begin{equation}\label{eq: con}
K(\ell)=\max\{ \kappa(G) \; ; \; \; |E(G)|=\ell\}
\end{equation}
\par
\noindent
and 

\begin{equation}\label{eq: conn}
k(\ell)=\min\{ \kappa(G) \; ; \; \; |E(G)|=\ell\}
\end{equation}
\par
\noindent
as well as 

\begin{equation}\label{eq: def}
K(v,\ell)=\max\{ \kappa(G) \; ; \; |V(G)|=v, \; |E(G)|=\ell\}
\end{equation}
 
We conjecture

\begin{conj}\label{conj:con}
For every $\ell\in \N$ 
$$K(\ell)=3^\ell.$$
\end{conj}

In fact, we will show that $K(\ell)$ lies between $3^\ell$ and $4^\ell$ for
every natural number $\ell$.
We will also see that $k(\ell)$ is linear in $\ell$. 

As we will explain, the values of $\kappa(P_r)$, where $P_r$ is the $r$-vertex
path, are relevant when
investigating colliding permutations. Giving a lower bound on $\kappa(P_4)$
the lower bound of (\ref{eq:prev}) will be improved to
$10^{n/4-O(1)}$.  
 
Also, we will discuss the following conjecture and its relation to
Conjecture~\ref{conj:clak}: 

\begin{conj}\label{conj:coco}
For every even $v\in \N$ 
$$K(v, v-1)={v+1 \choose {\left \lfloor \frac{v+1}{2} \right \rfloor}}.$$
\end{conj}

The concept of graph--different sequences from a fixed alphabet goes back to
Shannon's  
classical paper on zero--error capacity \cite{Sh}. This fundamental work has
inspired much of   
information theory ever since while in combinatorics it led Claude Berge to
define the intriguing class of perfect graphs, see
\cite{Ber}, cf. also \cite{BeRA}. As the reader knows, Berge's conjectures
about 
the structure of perfect graphs (cf. \cite{BeRA}) have had a tremendous impact
on the evolution of  
combinatorics and are by now important  and deep theorems at the center of the 
field. 
As explained in the survey  \cite{KO}, a large body of problems in extremal
combinatorics can be treated as zero--error    
problems in information theory. For the relationship of the present problems
to zero--error information theory we refer to \cite{KM}.  

\section{Natural graphs and infinite permutations}

Let $G$ be a natural graph and let again $\kappa(G)$ be the maximum
cardinality of a set of infinite permutations any two elements of which are
$G$--different, provided that this number is finite. It is easy to see that
this is always the case. Let $\chi(G)$ denote the chromatic number of graph
$G$. 

\begin{lem}\label{lem:fin}
For every natural graph 
$$\kappa(G) \leq (\chi(G))^{|V(G)|} $$
holds.
\end{lem}

\proof
Let us consider a proper coloring $c: V(G)\to \{1,\dots,\chi(G)\}$ of $G$. 
Let us write $v=|V(G)|$ and denote  by $W=[\chi(G)]^{\N}$ the set of infinite
sequences over the alphabet $\{1,\dots,\chi(G)\}$. Let $\mu$ be the uniform
probability measure  on $W$.
Let us consider a set $C$ of pairwise $G$--different permutations. 
We assign to any $\pi \in C$ the set $W(\pi)$ of all those
sequences of $W$ that for all $u\in V(G)$ have the element $c(u)$ in the
position where $\pi$ contains $u$. By our hypothesis on $C$, the sets $W(\pi)$
are pairwise disjoint for the different elements $\pi \in C$ whence, 
$$1=\mu(W)\geq \sum_{\pi \in C} \mu(W(\pi))=\sum_{\pi \in C}\chi(G)^{-v}=|C|
 \chi(G)^{-v}.$$ 
\hfill $\Box$

In the rest of this section we first investigate $K(\ell)$ and
$k(\ell)$. Subsequently our new lower bound on $\rho(n)$ will be proved.   

Let us denote by $S(G)$ the set of non--isolated vertices of the graph $G$.
We introduce a graph transformation that increases the value of $\kappa$.

\begin{prop}\label{prop:mon}

Let $F$ and $G$ be two graphs with $G$ obtained from $F$ upon deleting an
arbitrary edge in $E(F)$ followed by the addition of two new vertices to
$V(F)$ so that the latter form an additional edge in $G$. Then 
$$ \kappa(F )\leq \kappa(G).$$
\end{prop} 
\proof
Let us consider the $m=\kappa(F )$ pairwise $F$-different infinite permutations of an
arbitrary optimal configuration for $F.$ Let $t$ be large enough  
for the initial prefixes of length $t$ of these infinite sequences to be
pairwise $F$--different and let $q$ be the largest integer  
appearing in their coordinates. By the finiteness of $\kappa(F)$ such $t$ and
$q$ exist. Without restricting generality we can suppose
that the new edge of $G$ is  
$\{c, d\}$ with both $c$ and $d$ being strictly larger than $q$. 
We also suppose that the edge we will delete is $\{a,b\}\in E(F).$ Let us
now suffix to each of our  
sequences a new sequence of the same length $t$ where the suffix  to a
sequence $x_1x_2 \dots x_t$ is obtained from it by  
substituting every $a$ with  $c$ and every $b$ with  $d$ while the remaining
coordinates are defined in an arbitrary  
manner but in a way that the coordinates of the overall sequence of length
$2t$ be all different. Clearly, the $m$  
new sequences of length $2t$ are $G$--different. The rest is obvious, since we
can complete the new sequences to  
yield infinite permutations any way we like.
\hfill $\Box$

A straightforward consequence of the previous proposition is the following.

\begin{cor} \label{cor:mon}
$K(\ell)=\kappa(\ell K_2).$
\hfill$\Box$
\end{cor}

Thus we know that $K(\ell)$ is achieved by $\ell$ independent edges. 
It seems equally interesting to determine which graphs achieve $k(\ell).$
At first glance one might think that $S(F)\subseteq S(G)$ implies $\kappa(F)
\leq \kappa(G)$, but this is false.  
In particular, complete graphs do not have minimum $\kappa$ among graphs with
the same number of edges. Yet, determining their $\kappa$ value seems an
interesting problem. As we will see, the right guess for what graphs achieve
$k(\ell)$ turn out to be stars, at least for $\ell$ not too small. 
%Q: Is $\kappa$ monotone in $\nu$ or $\tau$? 
Below we will study the value of $\kappa$ for complete graphs,
stars, and paths. In particular, path graphs will take us back  
to the original puzzle about colliding permutations. 

\begin{prop}\label{prop:com}

For the complete graph $K_n$ on $n$ vertices

$${(n+1)!  \over 2} \leq \kappa(K_n ).$$

\end{prop} 
\proof
Consider the set of even permutations of $[n+1]$ and suppose $V(K_n)=[n].$ One
can observe 
that these permutations  are $K_n$--different. Indeed, if two arbitrary
permutations of $[n+1]$ are not  
$K_n$--different, then they differ only in positions in which for some fixed
$i \in [n]$ one has $n+1$ and the other  
has $i.$ Thus any of these two permutations can be obtained from the other by
exchanging the positions of   $n+1$ and the corresponding $i.$ But then the two permutations have different parity and in
particular they cannot both be even.  
In particular, the undesired relation does not occur between even permutations and this gives us ${(n+1)!  \over 2}$ permutations of $[n+1]$ that are
$K_n$--different.  
Next extend each of these permutations to infinite ones by 
suffixing the remaining natural numbers in an arbitrary order. 
\hfill $\Box$

\begin{prop}\label{prop:matc}

For the graph of $\ell$ independent edges we have
$$3^\ell\leq \kappa(\ell K_2)\leq 4^\ell.$$
\end{prop} 

\proof
Notice that the graph $\ell K_2$ has chromatic number two and its number of
vertices is $2\ell$, whence our upper bound follows by 
Lemma~\ref{lem:fin}. 

To prove
the lower bound, let us denote the edge set of our graph by 
$E(\ell K_2)=\{\{1,2\}, \{3,4\}, \dots \{2\ell-1,2\ell\}\}$. 
Consider the set of cyclic permutations $C_1=\{(1 2 \star), (2\star 1), (\star
1 2)\}$ and for every $1<i \leq \ell$ the sets  
$C_i$ obtained from $C_1$ by replacing $1$ with $2i-1$ and $2$ with $2i.$ It
is clear that for every $i\in [\ell]$ any two of the three ministrings in $C_i$
"differ" in the edge $\{2i-1, 2i\}$ of our graph $\ell K_2,$ meaning that they
have somewhere in the same position the two different endpoints  of this
edge. But this means that the $3^\ell$ strings in their cartesian product 
$$C=\times_{i=1}^\ell C_i$$
are pairwise $\ell K_2$--different as requested. Replacing the symbol $\star$
in 
our strings in an arbitrary
order with the different numbers from
$[3\ell]-[2\ell]$ we obtain $3^\ell$ permutations of $[3\ell]$ that
continue to be 
pairwise $\ell K_2$--different. The extension to infinite permutations is  
as always.
\hfill $\Box$

\medskip

The only infinite class of graphs for which we are able to completely determine
$\kappa$ are stars, i.e., the complete bipartite  
graphs $K_{1,r}$. We have

\begin{prop}\label{prop:star}

For every $r$
$$\kappa(K_{1,r})=2r+1.$$

\end{prop} 
\proof
By Lemma \ref{lem:fin} we know that $\kappa(K_{1,r})<\infty$. Let us denote
its value by $m$.  
Let us consider the vertices of $K_{1,r}$ to be the elements of $[r+1]$ and
let 1 be the "central"  
vertex of degree $r$.
It is obvious that in a set of $m$ sequences (infinite permutations) achieving
the maximum we are looking for    
all the sequences must have the central vertex 1 in a different position. Let
us consider our $m$ sequences as vertices of a directed graph $T$ in which
$(a,b)\in E(T)$ if  
the sequence corresponding to $a$ has a $j\in\{2,\dots,r+1\}$
in the same position
where the $1$ of  
the sequence corresponding to $b$ is placed. Then, by definition the directed
graph $T$ must  
contain a tournament, implying that
$$|E(T)|\geq {m \choose 2}.$$
On the other hand, every $a\in V(T)$ has at most $r$ outgoing edges. This means
that 
$$|E(T)|\leq mr.$$
Comparing the last two inequalities we get
$$m\leq 2r+1.$$
To prove a matching lower bound, consider the following set of permutations of
$[2r+1].$ 
For every $i\in [2r+1]$ let us define the coordinates of the $i^{\rm th}$
sequence 
$x_1(i)x_2(i)\dots x_{2r+1}(i)$ by  
$x_i(i)=1$ and, in general,  
$x_{i+j}(i)=j+1$ for any $0\leq j\leq r$ where all the coordinate indices
are considered modulo 
$2r+1$. The remaining coordinates are defined in an arbitrary manner so that
the resulting sequences define permutations of $[2r+1]$. It is easily seen
that this is a  
valid construction. In fact, observe that for any of our sequences the
"useful"  
symbols, those of $[r+1]$, corresponding to the vertices of the star graph,
occupy $r+1$  
"cyclically" consecutive coordinates, forming cyclical intervals. Since
$2(r+1)>2r+1$, these intervals are pairwise intersecting, and thus for any two
of them there must be a coordinate in the intersection for which the "left
end" of one of the intervals is contained in the other.  
The resulting permutations can be  
considered as prefixes of infinite permutations in the usual obvious way.
\hfill $\Box$

\medskip
Now we are ready to return to the problem of determining $k(\ell)$, at least
for large enough $\ell$.
The following easy lemma will be needed.

\begin{lem} \label{lem:szorzat}
If a finite graph $F$ contains vertex disjoint subgraphs $F_1,\dots,F_s$, then
$$\kappa(F)\ge \prod_{i=1}^s \kappa(F_i).$$
\end{lem} 

\proof
The proof is a straightforward generalization of the construction given in 
the proof of Proposition~\ref{prop:matc}.
Let $\hat C_i$ be a set of $\kappa(F_i)$ infinite sequences that are obtained
from $\kappa(F_i)$ pairwise $F_i$-different permutations of $\mathbb N$ by
substituting all natural numbers $i\notin V(F_i)$ by a $\star.$ As the
sequences in $\hat C_i$ contain only a finite number of elements different
from $\star$, we can take some finite initial segment of all these sequences
that already contain all elements of $V(F_i)$. Let $C_i$ be the set of these
finite sequences. Now consider the set $$C:=\times_{i=1}^s C_i$$ of finite
sequences that each contain all vertices in $\bigcup_{i=1}^s V(F_i)$ exactly
once and finitely many $\star$'s. These sequences are also pairwise
$F_i$-different for some $F_i$, thus they are pairwise
$F$-different. By construction, their number is $\prod_{i=1}^s \kappa(F_i)$. 
Extending them to infinite permutations of $\mathbb N$ the
statement is proved. 
\hfill$\Box$
\smallskip

It is straightforward from the previous lemma that if $F+G$ denotes the vertex
disjoint union of graphs $F$ and $G$ then $\kappa(F+G)\ge\kappa(F)\kappa(G)$. 
We do not know any example for strict inequality here. If equality was always
true that would immediately imply Conjecture~\ref{conj:con}.

Now we use Lemma~\ref{lem:szorzat} to prove our main result on $k(\ell)$.

\begin{prop} \label{min}
Let $G$ be a natural graph with $n:=|S(G)|>20$ and $|E(G)|=\ell$. Then
$$\kappa(G)\ge 2\ell+1.$$
The value of $k(\ell)$ is achieved by the graph $K_{1,\ell}$ whenever
$\ell>150$. 
\end{prop}

\proof
Let $G$ be a graph as in the statement 
and let $\nu=\nu(G)$ denote the size of a largest matching in
$G$. 

First assume that $\nu\ge n/4$. Then by Proposition~\ref{prop:matc} and  
the obvious monotonicity of $\kappa$ we have $$\kappa(G)\ge \kappa(\nu K_2)\ge
3^{\nu}.$$
Since $G$ is simple, we have $\ell\leq {n\choose 2},$ thus $k(\ell)\leq
\kappa\left(K_{1,{n\choose 2}}\right)=n(n-1)+1.$ So in this case (when $\nu\ge
n/4$) it is 
enough to prove that $$3^{\lceil n/4\rceil}\ge n(n-1)+1$$ holds. This is 
true if $n>20$.  

Next assume that $3\leq\nu<n/4$. Consider a largest matching of $G$ consisting
of edges $\{u_{2i-1},u_{2i}\}$ with $i=1,\dots,\nu$. The set $U:=\{u_1,\dots,
u_{2\nu}\}$ covers all edges of $G$ thus $\ell\leq {{2\nu}\choose
  2}+2\nu(n-2\nu).$ So we have $$k(\ell)\leq
\kappa(K_{1,\ell})\leq 2\left[{{2\nu}\choose 2}+2\nu(n-2\nu)\right]+1$$
in this case. On the other hand, for each vertex $a\in S(G)\setminus U$ there
is an edge $\{a,u_i\}$ for some $i$. We also know that if $a$ and $b$ are two
distinct vertices in $S(G)\setminus U$ and one of them is connected to
$u_{2j-1}$ (resp. $u_{2j}$) for some $j$, then the other one cannot be
connected to $u_{2j}$ (resp. $u_{2j-1}$) since 
otherwise replacing the matching
edge $\{u_{2j-1},u_{2j}\}$ with the other two edges of the path formed by the
vertices $a, u_{2j-1}, u_{2j}, b$ would result in a larger matching, a
contradiction. Choosing an edge for each $a\in
S(G)\setminus U$ that connects it to a vertex in $U$, we can form vertex
disjoint star subgraphs $K_{1,\ell_1},\dots,K_{1,\ell_{\nu}}$ of $G$, where
$\ell_i\ge 1$ for all $i$ and $\sum_{i=1}^{\nu}\ell_i=n-\nu.$ Then by
Lemma~\ref{lem:szorzat} and Proposition~\ref{prop:star} we have $\kappa(G)\ge
\prod_{i=1}^{\nu} (2\ell_i+1).$ The latter product is minimal (with respect to
the conditions on the $\ell_i$'s) if one $\ell_i$, say $\ell_1$,
equals to $n-2\nu+1$ and $\ell_2=\dots=\ell_{\nu}=1.$  Thus it is enough to
prove that $$3^{\nu-1}(2(n-2\nu+1)+1)> 2\left[{{2\nu}\choose
    2}+2\nu(n-2\nu)\right]+1$$  as the left hand side is a lower bound on
$\kappa(G)$ while the right hand side is an upper bound on $k(\ell)$. 
The latter inequality 
would be implied by
$$3^{\nu-1}(n-2\nu+1)>{{2\nu}\choose 2}+2\nu(n-2\nu)=\nu(2n-2\nu-1)$$ 
which, in turn, is equivalent to 
$${3^{\nu-1}\over\nu}>{{2n-2\nu-1}\over {n-2\nu+1}}$$
The left hand side of this last inequality is at least $3$ if $\nu\ge 3$,
while the right hand side is strictly less than $3$ for $\nu\leq n/4.$

The only case not yet covered is that of $\nu<3$. For $\nu=1$ there is nothing
to prove since then $G$
 itself is a star.
 If $\nu=2,$ then let a largest
matching be formed by the two edges $\{u_1,u_2\}$ and $\{u_3,u_4\}$, while 
once  again let
$U$ 
denote the union of their vertices. 
Let $a_1,\dots,a_{n-4}$ be the rest of the non-isolated vertices of $G$ and note that
$n-4>16$. Assume some
$a_i$ is connected to both $u_1$ and $u_2$ yielding a triangle. 
Then no $a_j, j\neq i$ can be
connected to either of $u_1$ or $u_2$, otherwise we could form a larger
matching. For similar reasons, if any $a_j$ is connected to $u_3$ then no 
$a_s, s\neq j$
can be connected to $u_4$. (If some $a_i$ forms a triangle with $u_1, u_2$ and
some $a_j$ with $u_3$ and $u_4$, then the remaining 
vertices $a_s$ 
must be isolated implying $n\leq 6$, a contradiction.) 
Thus if $a_i$ is connected to both $u_1$ and $u_2$, then the rest of the
$a_j$'s form a star centered at either $u_3$ or $u_4$. Thus in this case,
using  
again Lemma~\ref{lem:szorzat}, Propositions~\ref{prop:com} and \ref{prop:star}
imply 
$\kappa(G)\ge 12[2(n-4)+1]=24n-84$. The foregoing also implies $\ell\leq n+4$,
thus $k(\ell)\leq 2n+9<24n-84,$ whenever $n\ge 5$. Clearly, the situation
is similar if we exchange the role of the two matching edges. 

Assuming that no triangle is formed, we can again attach each vertex in
$S(G)\setminus U$ to one of the edges $\{u_1,u_2\}$ and $\{u_3,u_4\}$,
whichever it is connected to. Two vertex disjoint stars can be formed this way
establishing the lower bound $\kappa(G)\ge 3(2(n-3)+1)=6n-15$. For the number
of edges we 
now get $\ell\leq 6+2(n-4)=2n-2$ 
since the graph induces at most 6 edges on $U.$  
Thus we have
$k(\ell)\leq 4n-3$ which is less then $6n-15$ if $n>6$. 
This completes the proof of the first statement. 

\smallskip
If a simple graph has at most $20$ vertices then its number of edges is at
most $190$, so the second statement immediately follows from the first one if
$\ell>190$. If the graph contains a $K_6$ subgraph, then by
Proposition~\ref{prop:com} we have $\kappa(G)\ge 7!/2>381=\kappa(K_{1,190})\ge
k(\ell)$ if $\ell\leq 190$. Thus we may assume $K_6\nsubseteq G$ and this implies  by
Tur\'an's theorem that $\ell\leq 160$ if $n\leq 20$. But
$\kappa(K_{1,160})=321\leq 6!/2$, so if the conclusion is not true, we may
also assume that $G$ has no $K_5$ subgraph. Applying Tur\'an's theorem again,
this gives $\ell\leq 150$ for $n\leq 20$. Thus the statement is true whenever
$\ell>150$.  
\hfill$\Box$

\smallskip

\par\noindent{\bf Remark 1.}
We are quite convinced that the statement of Proposition~\ref{min} holds
without any restriction on $n$ or $\ell$. Some improvement on our treshold on
$\ell$ is easy to obtain. It seems to us, however, that
proving the statement in full generality either leads to tedious case
checkings or needs some new ideas.  
\hfill$\Diamond$

\medskip

The problem of determining $\kappa$ seems interesting in itself, moreover, it
helps to obtain better bounds for  
the original question on colliding permutations.  To explain this, 
we introduce a notion connecting the two questions. Let  
$\kappa(G,n)$ be the maximum number of pairwise $G$-different
permutations of $[n]$. Clearly, 

\begin{equation}\label{eq: bri}
\kappa(G)=\sup_{n} \kappa(G,n).
\end{equation}
Notice that by the finiteness of $\kappa(G)$ the supremum above is always
attained, so we could write maximum instead. 
Further, for the  graph $P_r$, the path on $r$ vertices, we have the
following.  
\begin{lem}\label{lem: pa}
For every $n>m>r$ the function $\rho$ satisfies the recursion

$$\rho(n)\geq \kappa(P_r, m) \rho(n-r).$$
\end{lem} 

\proof
We will call two arbitrary sequences of integers colliding if they have the same length
and if 
somewhere in the  
same position they feature integers differing by $1.$   
By the definition of $\kappa(P_r, m)$ we can construct this many sequences of
length $m$ such that in each of them every vertex of  
$P_r$ appears exactly once, the other positions are occupied by the "dummy"
symbol $\star$ and moreover these sequences are pairwise $P_r$--different.  
The latter implies that these sequences are pairwise colliding.
Furthermore, we have, also by definition, $\rho(n-r)$ permutations of $[n-r]$
that are pairwise colliding. Let us "shift" these permutations  
by adding $r$ to all of their coordinates. The new set of permutations of the
set $r+[n-r]=[r+1, n]$ maintains the property that its elements  
are pairwise colliding. Next we execute our basic operation of 
``substituting'' the
permutations of the second set into those coordinates of any sequence  
${\mbf x}$ from the first set where the sequence ${\mbf x}$ has a star. More
precisely, consider any sequence ${\mbf x}=x_1x_2 \dots x_m$ from our first
set and let  
$S({\mbf x})\in {[m] \choose m-r}$ be the set of those coordinates which are
occupied by stars. Let further ${\mbf y}=y_1y_2 \dots y_{n-r}$ be an  
arbitrary sequence from our second set, i.e., a permutation of $[r+1,n].$ The
sequence ${\mbf z}={\mbf y}\rightarrow {\mbf x}$ is a sequence of length $n$ in
which the first $m$ coordinates are defined in the following manner. We have
the equality $z_i=x_i$, if $i\leq m$ and $i\not \in S({\mbf x})$.  
Suppose further that $S({\mbf x})=\{j_1, j_2, \dots, j_{m-r}\}$. In the
$j_k^{\rm th}$ 
position we replace the symbol $\star$ by $y_k.$ (For $i>m$ we 
set
$z_i=y_{i-r}$.) Clearly, the resulting  
sequence is a permutation of $[n].$ Further, the so obtained $\kappa(P_r, m)
\rho(n-r)$ permutations are pairwise colliding. 
\hfill $\Box$

Observe next the following equality.

\begin{lem}\label{lem:puz}
$$\kappa (P_4,5)=10.$$
\end{lem} 

\par\noindent{\bf Remark 2.}
The existence of 10 permutations of $\{1\dots,5\}$ with the requested
properties is implicit in  
\cite{KM} since the construction of the 35 colliding permutations of 
$\{1,\dots,7\}$ in that paper does contain such a set in some appropriate
projection of its coordinates. 
Still, we prefer to give a simpler direct proof here. The argument for the
reverse inequality is the same as that proving $\rho(n)\leq {n\choose
  {\lfloor{n/2}\rfloor}}$ in \cite{KM}. 
\hfill$\Diamond$
\medskip

\proof
Let us consider the 10 permutations of $\{1,\dots,5\}$ obtainable by
considering the cyclic configurations of $(1,2)$, $(4,3)$ and the single
element $5$.  
We indeed have 10 different permutations by "cutting" in all the 5 possible 
ways both of the two cyclic configurations 3 building blocks can define. 
(So these are $12435, 24351, 43512, 35124, 51243$ and similarly the five
cyclic shifts of the sequence $43125$.) 
Let us further consider the graph $P_4$ (or, in fact, $P_4+K_1$) with vertex
set $\{1,\dots,5\}$ and with edge 
set $\{\{1,2\}, \{2,3\}, \{3,4\}\}$. In other words, consecutive numbers are  
adjacent vertices but $5$ is isolated. It is easy to check that the 
10 sequences above are $P_4$--different for the natural graph we defined. 
One can verify this by hand, 
yet let us give a more structured argument. 

The statement is true for two of our permutations if 
any of the two blocks with two elements are featured in intersecting positions
in 
the two respective permutations, or if  the two different blocks of length two  
are completely overlapping at least once. We claim  that one of these two
things will always happen. In  
fact, suppose  to the contrary to  have two different cyclic configurations,
say red and blue so that the only possible intersections  
of their blocks are intersections in one element between a  $(1,2)$
and a $(4,3)$ of different colors. This implies that  
the four cyclic intervals of length two are contained in a cycle of length 5
with only two points covered twice. But this is  
impossible as their total length is 8 and 8-5 is strictly larger than 2.

To see that $10$ is an upper bound it is enough to observe that the two even
elements of $\{1,\dots,5\}$ cannot be placed in the same two positions in two
permutations belonging to a set of $P_4$-different permutations of
$\{1,\dots,5\}$.  
\hfill $\Box$

\medskip
The above construction gives the following improved lower bound for the
exponential asymptotics of $\rho(n).$  

\begin{prop}\label{prop:as}

$$\lim_{n \rightarrow \infty} \rho^{1 \over n}(n) \geq 10^{1\over 4}$$
\end{prop} 
\proof
A simple combination of our two preceding lemmas implies 
$$\rho(n)\geq 10 \rho(n-4).$$
An iterated application of this inequality gives the desired result.

\hfill $\Box$

To close this section, let us take another look at Lemma \ref{lem:puz}. 
We believe that in fact $\kappa(P_4)=\kappa (P_4,5)=10$ 
and more generally, 
$$\kappa(P_v)=\kappa(P_v, v+1)={v+1 \choose {\left \lfloor \frac{v+1}{2}
    \right \rfloor}}$$ 
for even values of $v$. 
The original conjecture (see Conjecture~\ref{conj:clak}) for $\rho(v)$ would
be an immediate consequence of this  conjecture. 
To see this, suppose first that $v$ is even. Then
$\rho(v+1)=\kappa(P_{v+1},v+1)\geq \kappa(P_v, v+1)$
and  this would imply Conjecture~\ref{conj:clak} for odd values of $n$ right
away.  
Now, since for even $n$
$$\rho(n)= \kappa(P_n,n) \geq 2\kappa(P_{n-2}, n-1)$$
and likewise,
$${n \choose {n \over 2}}=2{{n-1} \choose {\left \lfloor \frac{n-1}{2} \right
    \rfloor}}$$ 
the last two relations would lead us to settle the conjecture for $n$
even. (The inequality above follows by putting an $n$ to the  end
of each sequence in an optimal construction for $\kappa(P_{n-2},n-1)$ and then
double each sequence by considering also its variant one obtains by exchanging
in it $n-1$ and $n$.) 

We also believe that
$$K(v,v-1)=\kappa(P_v).$$
As a combination of the two conjectures above we arrive at Conjecture
\ref{conj:coco} 
mentioned in the introduction.

\section{Related problems}

\subsection{A graph covering problem}

We show that 
determining $K(\ell)$
is equivalent to a graph covering problem introduced below.
The following standard definition is needed. 

\begin{defi} 
The (undirected) line graph $L(D)$ of the directed graph $D=(V,A)$ is defined by 
\begin{eqnarray*}
V(L(D))&=&A\\
E(L(D))&=&\{\{(a,b),(c,d)\}: b=c\ {\rm or}\ a=d\}
\end{eqnarray*}
\end{defi} 

Let ${\cal L}$ denote the family of all finite simple graphs that are
isomorphic to the line graph of some directed graph with possibly multiple
edges. It is a standard 
combinatorial problem 
to ask
how many graphs belonging to a certain
family of 
graphs 
are needed to cover all edges of a 
given 
complete graph, 
see, e.g., \cite{Khash}. We
show that the problem of determining $K(\ell)$ is equivalent to this problem
for the 
family ${\cal L}$. 

Let the minimum number of graphs in ${\cal L}$ the edge sets of which
 together can cover the edges of the complete graph $K_n$ be denoted by  
$h(n)$. 

\begin{prop} \label{lcover}
For any $M\in \mathbb N$, 
the minimum number $\ell$ for which $K(\ell)\ge M$ 
is equal to $h(M).$ 

\end{prop}  

\proof
Consider a construction attaining $K(\ell)$, that is a graph $G$ with $\ell$
edges and $K(\ell)$ infinite permutations that are $G$-different. Let this set
of permutations be denoted by $W$ and let $\{a,b\}$ be 
one of the edges of $G$. Define a graph $T_{a-b}$ on $W$ as its vertex set
where an edge is put between two permutations if and only
if there is a position where one of them has $a$ while the
other has $b$. In other words, the two permutations are $G$-different by the
edge $\{a,b\}$. 
Consider the graphs $T_{a-b}$ for all edges of $G$. These 
all have the same vertex set, while the union of their edge sets clearly
covers the complete graph $K_M$, where $M=K(\ell)$. 

Next we show that all the graphs 
$T_{a-b}$ belong to ${\cal L}$. 
To this end fix an edge $\{a,b\}\in V(G)$ and consider a graph $D_{a-b}$ with
its vertex set $V(D_{a-b})$ consisting of those positions where any of the
permutations 
in $W$ has a non-isolated vertex of $G$. Since $G$ and $W$ are finite,
so is $V(D_{a-b})$.
For each element of $W$ we define an edge of $D_{a-b}$. For
$\sigma\in W$, let $i$ and $j$ be the two positions where $\sigma$ contains
$a$ and $b$, respectively. Then let $\sigma$ be represented by the directed
edge $(i,j)$ in $D_{a-b}$. (If there is another permutation in $W$ with $a$
and $b$ being in the same positions as in $\sigma$ then we have another arc
$(i,j)$ in $D_{a-b}$ for this other permutation. Thus $D_{a-b}$ is a directed
multigraph.) Now it follows directly from the definitions that
$T_{a-b}=L(D_{a-b})$, thus $T_{a-b}$ is indeed the line graph of a digraph.  
Together with the previous paragraph this proves $h(K(\ell))\leq \ell.$
\medskip

For the reverse inequality consider a covering of $K_M$ with $h(M)$ graphs
belonging to ${\cal L}$. Let the line graphs in this covering be
$L_1,\dots,L_{h(M)}.$ We may assume that $V(L_i)=[M]$ for all $i$ by
extending the smaller vertex sets  
through the addition of 
isolated points. Let
$D_1,\dots,D_{h(M)}$ be 
directed graphs satisfying $L_i=L(D_i)$ for all $i$. (Such $D_i$'s exist since 
$L_i\in {\cal L}$.) By $E(D_i)=V(L_i)=[M]$ we can consider the edges of all
$D_i$'s labelled by $|E(D_i)|$ elements of $1,\dots,M$. (If $L_i$ had some
isolated vertices then the corresponding labels are not used.) 
Using these digraphs we define $M$ permutations $\sigma_1,\dots,\sigma_M$ 
that are $G$-different for the graph $G=\ell K_2$ with $\ell=h(M)$. For all
$i$ define $t_i=|V(D_i)|$ and identify $V(D_i)$ with $[t_i]$. Consider
$D_1$. If $D_1$ has an edge labelled $r$ and this edge is $(i,j)$, then put a
$1$ 
in position $i$ of $\sigma_r$ and put a $2$ in position $j$ of
$\sigma_r$. Do 
similarly for all edges of $D_1$. Then consider $D_2$. If it has an edge
labelled $r$ which is $(i',j')$ then put a $1$ in position $t_1+i'$ of
$\sigma_r$ 
and put a $2$ in position $t_1+j'$ of $\sigma_r$. In general, if $D_s$ has an
edge labelled $r$ which is $(a,b)$ then put a $2s-1$ 
in position
$(\sum_{k=1}^{s-1} t_k)+a$ and a $2s$ in position $(\sum_{k=1}^{s-1} t_k)+b$
of $\sigma_r$. When this is done for all edges of all $D_i$'s then extend the
obtained partial sequences to infinite permutations of $\mathbb N$ in an
arbitrary manner. This way 
one obtains $M$ permutations that are pairwise
$G$-different. To see this consider two of these permutations, say, $\sigma_q$
and 
$\sigma_r$. Look at the edge $\{q,r\}$ of our graph $K_M$ that was covered by
line graphs. Let $L_i$ be the line graph that covered the edge $\{q,r\}$. Then
$D_i$ has an edge labelled $q$ and another 
one labelled $r$ in such a way that the
head of the one is the tail of the other. This common point of these two edges
defines a position of $\sigma_q$ and $\sigma_r$ where one of them has $2i-1$
while the other has $2i$ making them $G$-different. 
\hfill$\Box$

\medskip
\par\noindent{\bf Remark 3.}
We note that the first part of the above proof makes no reference to the graphs
$\ell K_2$, yet it leads to another proof of the inequality $K(\ell)\leq
4^{\ell}$. Our earlier proof of this fact in
Proposition~\ref{prop:matc} relied on Corollary~\ref{cor:mon}. Here we sketch
a different proof.  
By Proposition~\ref{lcover} it
is enough to prove $h(M)\ge \log_4 M.$ Consider a line graph $L$ of a digraph
$D$ with $|V(D)|=t, |E(D)|=M$. Let $\hat K_t$ be the directed graph on $t$
vertices having an edge between any two different vertices in both
directions. $D$ can certainly be obtained by deleting some (perhaps zero)
edges of $\hat K_t$ and multiplying some (perhaps zero) of its edges. Thus
$L(D)$ can be obtained by multiplying some vertices of a subgraph
of $L(\hat K_t)$. (Multiplying a vertex means substituting it by an
independent set of size larger than $1$ in such a way that the
out-neighbourhoods and in-neighbourhoods of all vertices in this independent
set 
are the same as the corresponding neighbourhood of the original vertex.) 
This implies that the fractional chromatic number $\chi_f(L(D))$ of $L(D)$ is
bounded from above by $\chi_f(L(\hat K_t))$. (For the notion and basic
properties of the fractional chromatic number we refer to \cite{SchU}.) The
graph $L(\hat K_t)$ is vertex transitive so its fractional chromatic number is
equal to $|V(L(\hat K_t))|/\alpha(L(\hat K_t)),$ where $\alpha(F)$ stands for
the independence number of graph $F$. The latter ratio is bounded from above
by $4$ as $\hat K_t$ contains $\lfloor t/2\rfloor\cdot\lceil t/2\rceil$ edges
that 
form a complete bipartite subgraph and give rise to pairwise independent
vertices in the line
graph. So $\chi_f(L(D))\leq 4$. Now let $L_1,\dots,L_h$ be a minimal
collection of line graphs (of directed graphs) covering $K_M$. It is easy to
show that we must have $\prod_{i=1}^h \chi_f(L_i)\ge \chi_f(\bigcup_{i=i}^h
L_i)\ge \chi_f(K_M)=M$. Having $\chi_f(L_i)\leq 4$ for all $i$ this implies
$h\ge \log_4 M.$ 
\hfill$\Diamond$

\subsection{Fixed suborders}

It seems worthwile to revisit the problem of the determination of $\kappa(G)$
for the restricted class of infinite permutations in which the vertices of $G$
appear in a predetermined order. We will study this problem in
the case of complete graphs. Without restricting generality, we can
suppose that the fixed order is the natural one. 

Let $\kappa_{id}(K_n)$ denote the 
maximum number of infinite
permutations of $\mathbb N$ that are $K_n$-different and contain the first $n$
positive integers in their natural order. 

\begin{prop} \label{catalan}
For every $n\in \mathbb N$ $$\kappa_{id}(K_n)\ge C_n$$
holds, where $C_n={1\over{n+1}}{2n\choose n}$ is the $n^{\it th}$ Catalan
number. 
\end{prop}  

\proof
For $n=1,2$ we have equality: $\kappa_{id}(K_1)=1,
\kappa_{id}(K_2)=2$. 
Set $a_0=1$ and $a_n:=\kappa_{id}(K_n)$.
It is enough to prove that the numbers $a_n$
satisfy the inequality 
$$a_{n+1}\ge \sum_{i=0}^n a_i a_{n-i}$$
that has the well-known recursion of Catalan
numbers on its right hand side. 
We will look at our infinite permutations as infinite sequences consisting of
infinitely many $\star$'s and one of each of the symbols
$1,2,\dots,n$, where the $\star$'s refer to all other symbols. Clearly, only
the positions of the elements of $[n]$ are relevant
with respect to the $K_n$-difference relation. Thus we will define the
positions of the elements of $[n]$ and then let the $\star$'s be substituted
by the other numbers in any way 
that will result in infinite permutations of
$\mathbb N$. 

Our construction is inductive. Assume that we already know that $a_k\ge C_k$
holds for $k\leq n$ and 
thus it suffices to prove it for $n+1$. 
Fix a position of our permutations which is ``far away'', meaning that it is
far enough for having enough earlier positions for the following
construction. Call this position $j$. For each $i=0,\dots, n$ we construct 
$a_ia_{n-i}$ sequences having $i+1$ at their position $j$. Any two of these
sequences that
have a different symbol at position $j$ are $K_n$-different. 
For those sequences that have $i+1$ at their position $j$ do the
following. Consider a construction of $a_i$ pairwise $K_i$-different 
sequences consisting of symbols $1,\dots,i,\star$, where the symbols in $[i]$
are all used somewhere in the first $j-1$ positions (this is possible if $j$
is chosen large enough). Take the first $j-1$ coordinates of all
these sequences, $a_{n-i}$ times each, and continue each of them with an
$i+1$ at the $j^{\rm th}$ position. So we have $a_ia_{n-i}$ sequences of length
$j$ with $i+1$ at the $j^{\rm th}$ position, each of these sequences are one of
$a_i$ possible types and we have $a_{n-i}$ copies from each type. 

Now consider $a_{n-i}$ sequences with the
symbols $1,\dots,n-i, \star$ that are pairwise  $K_{n-i}$-different and shift
each value in these sequences by $i+1$. (The latter means that we change each
value $k$ to $k+i+1$ in these sequences while $\star$'s remain $\star$'s.) 
For each type of the previous sequences take its $a_{n-i}$ copies and suffix to
each of them one of the current $a_{n-i}$ different sequences. This way
one gets $a_ia_{n-i}$ $K_{n+1}$-different sequences with symbol $i+1$ at
position $j$. Doing this for all $i=1,\dots,n$ one obtains $\sum_{i=0}^n a_i
a_{n-i}$ 
$K_{n+1}$-different sequences proving the desired inequality. 
\hfill$\Box$

\medskip
\par\noindent
Observe that if we have a construction of $M$ infinite permutations that are
pairwise $K_{n+1}$-different and 
furthermore each of them contains the symbols
$1,\dots,n$ in 
their natural order then the number of 
those among them that  
have the symbol $i+1$ in a fixed position is at most
$\kappa_{id}(K_i)\kappa_{id}(K_{n-i})$. This is simply because all such
sequences must be made $K_{n+1}$-different entirely either by their smallest
$i$ or by their largest $(n-i)$ ``non-dummy'' symbols. So in case we have a
construction where at some position every 
permutation has a useful value, that is a natural number at most $n+1$, then
the number of these 
permutations is at most $C_{n+1}$ and the construction in the proof of
Proposition~\ref{catalan} is optimal. It seems plausible 
that the condition on this special coordinate can be dropped, implying that
Catalan 
numbers give the true optimum. If it is so then one feels it should be
possible to find a bijection between our permutations and the objects of one
of the many enumeration problems leading to Catalan numbers,
cf. \cite{Stanl}. 
It seems to be a
significant difficulty, however, that our permutations are not objects having
some  
structural property on their own, 
as it happens in most of the
enumeration problems leading to Catalan numbers. Rather, in our case the
criterion is in terms of a relation between pairs of objects, and this seems to make
an important difference.

\section{Acknowledgment}

The first author would like to thank Alexandr Kostochka for a stimulating
discussion. Thanks also to Riccardo Silvestri for computer help.

\end{document}